\documentclass[10pt]{amsart}
\usepackage{latexsym,amssymb}

\theoremstyle{plain}
  \newtheorem{theorem}{Theorem}[section]
  \newtheorem{proposition}[theorem]{Proposition}
  \newtheorem{lemma}[theorem]{Lemma}
  \newtheorem{corollary}[theorem]{Corollary}
  
\theoremstyle{definition}
  \newtheorem{definition}[theorem]{Definition}
  \newtheorem{example}[theorem]{Example}
  
  \newtheorem{question}[theorem]{Question}
 \theoremstyle{remark}
  \newtheorem{remark}[theorem]{Remark}

\numberwithin{equation}{section}

\def\NN{{\mathbb N}}
\def\ZZ{{\mathbb Z}}
\def\QQ{{\mathbb Q}}
\def\RR{{\mathbb R}}

\def\OOO{{\mathcal O}}
\def\CC{{\mathbb C}}

\def\symm{\mathfrak{S}}

\def\maj{\mathrm{maj}}
\def\fmaj{\mathrm{fmaj}}
\def\des{\mathrm{des}}
\def\fdes{\mathrm{fdes}}
\def\inv{\mathrm{inv}}
\def\SYT{\mathrm{SYT}}
\def\Hilb{\mathrm{Hilb}}
\def\Res{\mathrm{Res}}
\def\Ind{\mathrm{Ind}}
\def\Des{\mathrm{Des}}

\def\Cent{\mathrm{Cent}}
\def\Gal{\mathrm{Gal}}
\def\diag{\Delta}

\def\maj{\mathrm{maj}}

\def\nneg{\mathrm{neg}}

\def\xx{\mathbf{x}}
\def\yy{\mathbf{y}}

\begin{document}

\title[Bimahonian distributions]
{Bimahonian distributions}

\author{H\'el\`ene Barcelo}
\email{barcelo@asu.edu}
\address{ Department of Mathematics and Statistics\\
Arizona State University\\
Tempe, AZ 85287}

\author{Victor Reiner}
\email{reiner@math.umn.edu}

\author{Dennis Stanton}
\email{stanton@math.umn.edu}
\address{School of Mathematics\\
University of  Minnesota\\
Minneapolis, MN 55455}

\thanks{First author supported by NSA grant H98230-05-1-0256. 
Second and third authors supported by NSF grants DMS-0601010 and DMS-0503660, respectively.}



\begin{abstract}
Motivated by permutation statistics, we
define for any complex reflection group $W$ a family of bivariate generating
functions $W^\sigma(t,q)$.
They are defined either in terms of Hilbert series for $W$-invariant polynomials when 
$W$ acts diagonally on two sets of variables, or equivalently, as sums
involving the fake degrees of irreducible representations for $W$.
It is shown that $W^\sigma(t,q)$ satisfies a ``bicyclic sieving phenomenon'' which 
combinatorially interprets its values when $t$ and $q$ are certain roots of unity.
%
\end{abstract}

\maketitle

\tableofcontents

\section{Introduction}
\label{intro-section}

  This paper contains one main new definition, of the {\it bimahonian distributions} for a complex 
reflection group (Definition~\ref{bimahonian-fake-degree-defn}), 
and three main results about it 
(Theorems~\ref{bimahonian-invariant-defn}, \ref{bimahonian-by-fmaj-theorem}, \ref{biCSP-theorem}),
inspired by known results in the theory of permutation statistics.  We explain here some background
and context for these results.

\subsection{The mahonian distribution}

P.A. MacMahon \cite{Macmahon} proved that there are 
two natural equidistributed statistics on permutations $w$ in 
the symmetric group $\symm_n$ on $n$ letters, namely
the {\it inversion number} (or equivalently, the Coxeter group length)
\begin{equation}
\label{inv-definition}
\inv(w):=|\{(i,j):1 \leq i < j \leq n, w(i)> w(j) \} \,\, (=\ell(w))
\end{equation}
and the {\it major index} 
\begin{equation}
\label{maj-definition}
\maj(w):=\sum_{i \in \Des(w)} i 
\end{equation}
where $\Des(w):=\{ i : w(i) > w(i+1) \}$ is the {\it descent set} of $w$.

Foata \cite{Foata} dubbed their common distribution 
\begin{equation}
\label{type-A-mahonian}
\sum_{w \in \symm_n} q^{\inv(w)}= \sum_{w \in \symm_n} q^{\maj(w)} = \prod_{i=1}^n (1+q+q^2+ \cdots + q^{i-1})
\end{equation}
the {\it mahonian distribution}.  Subsequently many other statistics on $\symm_n$
with this distribution have been discovered; see e.g., Clarke \cite{Clarke}.

This distribution has a well-known generalization to any {\it complex reflection
group} $W$, by means of invariant theory, as we next review.  We also review
here some known generalizations of the statistics $\inv(w), \maj(w)$,
which are unfortunately currently known in less generality.

Recall that a 
{\it complex reflection group} is a finite subgroup of $GL(V)$ for $V=\CC^n$ generated by {\it reflections}
(= elements whose fixed subspace has codimension $1$).  Shepard and 
Todd \cite{ShepardTodd} classified such groups.  They also
showed (as did Chevalley \cite{Chevalley}) that they are characterized among the
finite subgroups $W$ of $GL(V)$ as those whose action on the symmetric algebra 
$S(V^*)=\CC[V] \cong \CC[x_1,\ldots,x_n]$
has the invariant subalgebra $\CC[V]^W$ again a polynomial algebra.
In this case one can choose homogeneous generators $f_1,\ldots,f_n$ for $\CC[V]^W$,
having {\it degrees} $d_1,\ldots,d_n$, and define the {\it mahonian distribution} for $W$, in analogy with the product formula in \eqref{type-A-mahonian}, by
$$
W(q):=\prod_{i=1}^n (1+q+q^2+ \cdots + q^{d_i-1}).
$$
This polynomial $W(q)$ is a ``$q$-analogue'' of $|W|$ in the sense that
$W(1)=|W|$.  It is a polynomial in $q$ whose degree $N^*:=\sum_{i=1}^n (d_i-1)$
is the number of reflections in $W$.

We now give two equivalent ways to rephrase the definition.  
$W(q)$ is also the {\it Hilbert series}
for the {\it coinvariant algebra} $\CC[V]/(\CC[V]^W_+)$:
\begin{equation}
\label{coinvariant-mahonian-defn}
W(q)=\Hilb( \,\, \CC[V]/(\CC[V]^W_+) \,\, ,q)
\end{equation}
where $(\CC[V]^W_+):=(f_1,\ldots,f_n)$ is the ideal in $\CC[V]$ generated by the invariants $\CC[V]^W_+$
of positive degree.  In their work, both Shepard and Todd \cite{ShepardTodd} and
Chevalley \cite{Chevalley} showed that the coinvariant algebra $\CC[V]/(\CC[V]^W_+)$ carries the {\it regular
representation} of $W$.  Consequently, each irreducible complex character/representation
$\chi^\lambda$ of $W$ occurs in it with multiplicity equal to its
degree $f^\lambda:=\chi^\lambda(e)$.  This implies that
there is a unique polynomial $f^\lambda(q)$ in $q$ with
\begin{equation}
\label{STC-rephrased}
f^\lambda(1)=f^\lambda:=\chi^\lambda(e)
\end{equation}
called the {\it fake degree polynomial}, 
which records the homogeneous components of $\CC[V]/(\CC[V]^W_+)$
in which this irreducible occurs:
\begin{equation}
\label{fake-degree-definition}
f^\lambda(q) 
  := \sum_{ i \geq 0 } q^i
    \langle \,\, \chi^\lambda \,\, , \,\, \left( \CC[V]/(\CC[V]^W_+) \right)_i \,\, \rangle_W .
\end{equation}
Here $R_i$ denotes the $i^{th}$ homogeneous component of a graded ring $R$,
and $\langle \cdot ,\cdot \rangle_W$ denotes the inner product (or intertwining number) of 
two $W$-representations or characters.
Thus our second rephrasing of the definition of $W(q)$ is:

\begin{equation}
\label{fake-degree-mahonian-defn}
W(q)=\sum_\lambda f^\lambda(1) f^\lambda(q).
\end{equation}

At somewhat lower levels of generality, one finds formulae for $W(q)$ generalizing
the formulae involving the statistics $\inv(w), \maj(w)$ that appeared
in \eqref{type-A-mahonian}:
\begin{enumerate}
\item[$\bullet$]  When $W$ is a (not necessarily crystallographic
\footnote{If one further assumes that $W$ is crystallographic, and hence a Weyl group, 
the coinvariant algebra $\CC[V]/(\CC[V]^W_+)$ gives the cohomology ring of the 
associated {\it flag manifold} $G/B$, and $W(q)$ is also the {\it Poincar\'e series} for $G/B$.} 
) 
real reflection group, and hence a Coxeter
group, with Coxeter length function $\ell(w)$, one has
$$
W(q)=\sum_{w \in W} q^{\ell(w)}.
$$
\item[$\bullet$]  When $W=\ZZ/d\ZZ \wr \symm_n$, the wreath product of a
cyclic group of order $d$ with the symmetric group,  $W(q)$ is the distribution
for a statistic $\fmaj(w)$ defined by Adin and Roichman \cite{AdinRoichman}
(see Section~\ref{fmaj-section} below):
\begin{equation}
\label{fmaj-mahonian}
W(q) = \sum_{w \in W} q^{\fmaj(w)}.
\end{equation}
\end{enumerate}

\subsection{The bimahonian distribution}

Foata and Sch\"utzenberger \cite{FoataSchutzenberger} first observed that
there is a bivariate distribution on the symmetric group $\symm_n$, shared by
several pairs of mahonian statistics, including the pairs 
\begin{equation}
\label{bimahonian-by-maj-defn}
(\maj(w),\inv(w))\text{ and }(\maj(w),maj(w^{-1}))
\end{equation}
This bivariate distribution is closely related to the fake degrees for $W=\symm_n$
and {\it standard Young tableaux} via the Robinson-Schensted correspondence,
as well as to bipartite partitions and invariant theory;   
see Sections~\ref{Gordon-section} and \ref{fmaj-section} for more discussion and references.  

The goal of this paper is to generalize this ``bimahonian'' 
distribution on $W=\symm_n$ to any complex reflection group $W$.
In fact, one is naturally led to consider a family of such bivariate distributions
$W^\sigma(t,q)$, indexed by field automorphisms $\sigma$ lying in the Galois group
\cite[\S14.5]{DummitFoote}
$$
\Gal(\QQ[e^{\frac{2\pi i}{m}}]/\QQ) \cong (\ZZ/m\ZZ)^\times
$$ 
of any cyclotomic extension $\QQ[e^{\frac{2\pi i}{m}}]$ of $\QQ$ 
large enough to define the matrix entries for $W$.  It is known \cite[\S 12.3]{Serre} that one can
take $m$ to be the least common multiple of the orders of the elements $w$ in $W$. When referring
to Galois automorphisms $\sigma$, we will implicitly assume that $\sigma \in \Gal(\QQ[e^{\frac{2\pi i}{m}}]/\QQ)$
for this choice of $m$, unless specified otherwise.

\begin{definition}
\label{bimahonian-fake-degree-defn}
Given $\sigma$, in 
analogy to \eqref{fake-degree-mahonian-defn}, 
define the {\it $\sigma$-bimahonian distribution} for $W$ by
\begin{equation}
\label{actual-bimahonian-defn-eqn}
W^\sigma(t,q):=\sum_{\lambda} f^{\sigma(\lambda)}(t) f^{\overline{\lambda}}(q)
             =\sum_{\lambda} f^\lambda(t) f^{\overline{\sigma}^{-1}(\lambda)}(q).
\end{equation}
Here both sums run over all irreducible complex $W$-representations $\lambda$,
and 
$$
\chi^{\sigma(\lambda)}(w)=\chi^{\lambda}(\sigma(w))=\sigma(\chi^\lambda(w))
$$
denotes the character of the irreducible representation $\sigma(\lambda)$ 
defined by applying $\sigma$ entrywise to the matrices representing the group
elements in $\lambda$.  
\end{definition}

The two most important examples of $\sigma$ will be
\begin{enumerate}
\item[$\bullet$]
$\sigma(z)=z$, for which we denote $W^\sigma(t,q)$ by $W(t,q)$, and
\item[$\bullet$]
$\sigma(z)=\bar{z}$, where $\sigma(\lambda)$ is equivalent to the 
representation {\it contragredient} to $\lambda$, and for which 
we denote $W^\sigma(t,q)$ by  $\overline{W}(t,q)$.
\end{enumerate}

\noindent
From Definition~\ref{bimahonian-fake-degree-defn} 
one can see that 
\begin{equation}
\label{t-q-symmetry}
W^\sigma(q,t) = W^{\sigma^{-1}}(t,q)
\end{equation}
and hence both $W(t,q), \overline{W}(t,q)$ are symmetric polynomials in $t,q$.
Setting $t=1$ or $q=1$ in \eqref{actual-bimahonian-defn-eqn} and comparing with
\eqref{fake-degree-mahonian-defn} gives 
$$
W^\sigma(1,q)=W^\sigma(q,1)=W(q)
$$ 
for any $\sigma$.

  It turns out that, by analogy to \eqref{coinvariant-mahonian-defn},
one can also define $W^\sigma(t,q)$ 
as a {\it Hilbert series} arising in invariant theory.
One considers the {\it $\sigma$-diagonal} embedding of $W$ defined by
$$
\diag^\sigma W :=\{(w,\sigma(w)): w \in W\} \subset W \times W^\sigma \subset GL(V) \times GL(V)
$$
where $W^\sigma:=\{\sigma(w): w\in W\} \subset GL(V)$.
Note that $W \times W^\sigma$ acts on the symmetric algebra
$$
S(V^* \oplus V^*)=\CC[V \oplus V] \cong \CC[x_1,\ldots,x_n,y_1,\ldots,y_n]
$$
in a way that preserves the $\NN^2$-grading given by $\deg(x_i)=(1,0), \deg(y_j)=(0,1)$
for all $i,j$.  Hence one can consider the
$\NN^2$-graded Hilbert series in $t,q$ for various graded subalgebras and
subquotients of $\CC[V \oplus V]$.  The following 
analogue of \eqref{coinvariant-mahonian-defn} is proven in Section~\ref{invariants-section} below.

\begin{theorem}
\label{bimahonian-invariant-defn}
For any complex reflection group $W$ and Galois automorphism $\sigma$,
$W^\sigma(t,q)$ is the $\NN^2$-graded Hilbert series for the
ring 
$$
\CC[V \oplus V]^{\diag^\sigma W}/\left( \CC[V \oplus V]^{W \times W^\sigma}_+ \right).
$$
\end{theorem}

Note that when $W$ is a {\it real} reflection group, one has $\overline{W}=W$ and 
$\diag^\sigma(W)=W$, so that $\overline{W}(t,q)=W(t,q)$.  For example,
in the original motivating special case where $W=\symm_n$, the description of $W(t,q)$ 
via Theorem~\ref{bimahonian-invariant-defn} relates to
work on {\it bipartite partitions} by 
Carlitz \cite{Carlitz}, Wright \cite{Wright}, Gordon \cite{Gordon}, Roselle \cite{Ros}, 
Solomon \cite{Solomon}, and Garsia and Gessel \cite{GarsiaGessel}; 
see \cite[Example 5.3]{Stanley-invariants}.  

  Much of this was generalized from $W=\symm_n$ 
to the wreath products $W=\ZZ/d\ZZ \wr \symm_n$ in work
of Adin and Roichman \cite{AdinRoichman}  (see also Bagno and Biagioli \cite{BagnoBiagioli},
 Bergeron and Biagioli \cite{BergeronBiagioli}, Bergeron and Lamontagne \cite{BergeronLamontagne}, 
Chow and Gessel \cite{ChowGessel}, Foata and Han \cite{FoataHan}, and Shwartz \cite{Shwartz}).
We discuss some of this in Section~\ref{fmaj-section} below, and prove the following interpretation
for $W^\sigma(t,q)$, generalizing a result of Adin and Roichman \cite{AdinRoichman}.

\begin{theorem}
\label{bimahonian-by-fmaj-theorem}
For the wreath products $W=\ZZ/d\ZZ \wr \symm_n$ and
any Galois automorphism $\sigma \in \Gal(\QQ[e^{\frac{2\pi i}{d}}]/\QQ)$, one has
$$
W^\sigma(t,q) = \sum_{w \in W} q^{\fmaj(w)} t^{\fmaj(\sigma(w^{-1}))}.
$$
\end{theorem}

\subsection{Bicyclic sieving}

One motivation for the current work came from previous studies \cite{BMS, BSS} of
the coefficients $a_{k,\ell}(i,j)$ uniquely defined by
\begin{equation}
\label{reduced-coefficients-defn}
W^\sigma(t,q) \equiv \sum_{\substack{0 \leq i < k \\ 0 \leq j < \ell}} a_{k,\ell}(i,j) t^i q^j 
\mod (t^k-1, q^\ell-1)
\end{equation}
in the case where $W=\symm_n$ and $\sigma$ is the identity.
Equivalent information is provided by 
knowing the evaluations $W^\sigma(\omega,\omega')$ where as $\omega, \omega'$ vary over all
$k^{th}, \ell^{th}$ complex roots of unity, respectively; these were studied in
the case $W=\symm_n$ by Carlitz \cite{Carlitz} and Gordon \cite{Gordon}.  

In Section~\ref{Springer-section} below, 
we prove Theorem~\ref{biCSP-theorem}, generalizing some of these results to
all complex reflection groups, and providing a combinatorial interpretation for
some of these evaluations  $W^\sigma(\omega,\omega')$.  This is our first instance of a 
{\it bicyclic sieving phenomenon}, generalizing
the notion of a {\it cyclic sieving phenomenon} introduced in \cite{RSW}.

\begin{theorem}
\label{biCSP-theorem}
Let $C, C'$ be cyclic subgroups of $W$ generated by regular elements $c, c'$ in $W$,
having regular eigenvalues $\omega,\omega'$, in the sense of Springer \cite{Springer};
see \S\ref{Springer-section} below.
Given the Galois automorphism $\sigma$, choose an integer $s$ with the property that 
$\sigma(\omega)=\omega^s$, and then define a (left-)action of $C \times C'$ 
on $W$ as follows:
$$
(c,c')w := c^s w \sigma(c')^{-1}.
$$ 
Then for any integers $i,j$ one has 
$$
W^\sigma(\omega^{-i},(\omega')^{-j})=|\{w \in W: (c^i,(c')^j)w=w\}|.
$$
\end{theorem}

Alternatively, one can phrase this phenomenon (see Proposition~\ref{biCSP-defn-prop}) as a combinatorial
interpretation of the coefficients $a_{k,\ell}(i,j)$ in \eqref{reduced-coefficients-defn}:  
if the regular elements $c,c'$ in our complex reflection group $W$ have orders $k, \ell$, 
then $a_{k,\ell}(i,j)$ is the number of $C \times C'$-orbits on $W$ for which the  
stabilizer subgroup of any element in the orbit lies in the kernel of the character 
$\rho^{(i,j)}: C \times C' \rightarrow \CC^\times$ sending $(c, c') \mapsto \omega^{-i}(\omega')^{-j}$.


\section{Proof of Theorem~\ref{bimahonian-invariant-defn}}
\label{invariants-section}

After some preliminaries about Galois automorphisms $\sigma$, we recall 
the statement of Theorem~\ref{bimahonian-invariant-defn}, prove it,
and give some combinatorial consequences for the bimahonian distributions
$W^\sigma(t,q)$.

For any Galois automorphism 
$\sigma \in \Gal(\QQ[e^{\frac{2 \pi i}{m}}]/\QQ)$, the group $W^\sigma=\sigma(W)$
is also a complex reflection group:
$\sigma$ takes a reflection $r$ having characteristic polynomial
$(t-\omega)(t-1)^{n-1}$ to a reflection $\sigma(r)$ having
characteristic polynomial $(t-\sigma(\omega))(t-1)^{n-1}$.
Furthermore, if $\CC[V]^W=\CC[f_1,\ldots,f_n]$, 
then 
$$
\CC[V]^{W^\sigma}=\CC[\sigma(f_1),\ldots,\sigma(f_n)],
$$
and hence $W$ and $W^\sigma$ share the same degrees $d_1,\ldots,d_n$ for their basic invariants.
This implies
\begin{equation}
\label{conjugate-Hilbs-equal}
\Hilb(\CC[V]^W,q)=\Hilb(\CC[V]^{W^\sigma},q)=\prod_{i=1}^n \frac{1}{1-q^{d_i}}.
\end{equation}
This gives a relation between the graded character of 
the $W$-coinvariant algebra $A:=\CC[V]/(\CC[V]^W_+)$ and 
the $W^\sigma$-coinvariant algebra $A^\sigma:=\CC[V]/(\CC[V]^{W^\sigma}_+)$.
Let $A_j, A^\sigma_j$ denote their $j^{th}$-graded components.

\begin{proposition}
\label{coinvariant-trace-comparison}
For every $j$, one has 
$$
\chi_{A^\sigma_j}(\sigma(w)) = \sigma\left( \chi_{A_j}(w) \right)
$$
\end{proposition}
\begin{proof}
Since $\CC[V]^W, \CC[V]^{W^\sigma}$ are generated by homogeneous systems of parameters for
$\CC[V]$, and since $\CC[V]$ is a {\it Cohen-Macaulay} ring \cite[\S 3]{Stanley-invariants}, 
$\CC[V]$ is free as a module over either of these subalgebras.  Hence
$$
\begin{aligned}
\Hilb(\CC[V]^W,q) \left( \sum_j \chi_{A_j}(w) q^j \right) 
   & = \sum_j \chi_{\CC[V]_j}(w) q^j \\
\Hilb(\CC[V]^{W^\sigma},q) \left( \sum_j \chi_{A^\sigma_j}(\sigma(w)) q^j \right) 
   & = \sum_j \chi_{\CC[V]_j}(\sigma(w)) q^j \\
   & = \sigma \left( \sum_j  \chi_{\CC[V]_j}(w)  q^j \right) \\
\end{aligned}
$$
so that, dividing by $\Hilb(\CC[V]^W,q)(=\Hilb(\CC[V]^{W^\sigma},q))$, one obtains
$$
\begin{aligned}
\sum_j \chi_{A^\sigma_j}(\sigma(w)) q^j 
   & = \sigma \left( \frac{\sum_j \chi_{\CC[V]_j}(w) q^j}{\Hilb(\CC[V]^{W^\sigma},q)} \right) \\
   & = \sigma \left( \frac{\sum_j \chi_{\CC[V]_j}(w) q^j}{\Hilb(\CC[V]^{W},q)} \right) \\
   & = \sigma \left( \sum_j \chi_{A_j}(w) q^j \right) \\
   & =  \sum_j \sigma( \chi_{A_j}(w)) q^j.
\end{aligned}
$$
Here we have used throughout that the coefficients of any Hilbert series are integers,
and hence are fixed by $\sigma$,
while the second equality uses \eqref{conjugate-Hilbs-equal}.
Comparing coefficients of $q^j$ gives the assertion.
\end{proof}

We next recall the statement of Theorem~\ref{bimahonian-invariant-defn},
and prove it.

\vskip .2in
\noindent
{\bf Theorem~\ref{bimahonian-invariant-defn}.}
{\it 
For any complex reflection group $W$ and $\sigma \in \Gal(\QQ[e^{\frac{2 \pi i}{m}}]/\QQ)$,
$W^\sigma(t,q)$ is the $\NN^2$-graded Hilbert series for the
ring 
$$
\CC[V \oplus V]^{\diag^\sigma W}/\left( \CC[V \oplus V]^{W \times W^\sigma}_+ \right).
$$
}

\noindent
\begin{proof}
Note that the $W \times W^\sigma$-coinvariant algebra 
$$
A_{W \times W^\sigma}=\CC[V \oplus V]/(\CC[V \oplus V]^{W \times W^\sigma}_+)
$$ 
satisfies
$$ 
A_{W \times W^\sigma} \cong A \otimes A^\sigma
$$
as $W \times  W^\sigma$-representations, and hence for each $i, j$ one has
\begin{equation}
\label{tensor-decomposition}
\left( A_{W \times W^\sigma} \right)_{i,j} 
\cong 
A_i \otimes A^\sigma_j.
\end{equation}

Thus one has 
$$
\begin{aligned}
&\Hilb \left( \CC[V \oplus V]^{\diag^\sigma W} / (\CC[V \oplus V]^{W \times W^\sigma}_+); t,q \right) \\
&=\Hilb \left(  \left( A_{W \times W^\sigma} \right)^{\diag^\sigma W}; t,q \right) \\
&=\sum_{i,j} t^i q^j \left\langle
                              \left( A_{W \times W^\sigma} \right)_{i,j}, \mathbf{1}
                            \right\rangle_{\diag^\sigma W}.
\end{aligned}
$$
Here the first equality uses the fact that the {\it averaging (or Reynolds) operator}
$$
f \mapsto \frac{1}{|W|} \sum_{w \in W} w(f)
$$ 
gives an $\CC[V \oplus V]^{W \times W^\sigma}$-module projection
$\CC[V \oplus V] \twoheadrightarrow \CC[V \oplus V]^{\diag^\sigma W}$
that splits\footnote{We are being careful here-- the relation between the invariant
rings $k[V]^G, k[V]^H$ and $k[V]$ and coinvariants for subgroups $H \subset G \subset GL(V)$ 
is not as straightforward when working over a field $k$ of positive characteristic.
One always has a natural map $k[V]^H/(k[V]^G_+) \rightarrow \left( k[V]/(k[V]^G_+)\right)^H$,
but it can fail to be an isomorphism when one lacks a $k[V]^G$-module splitting 
(such as the Reynolds operator above) for the inclusion $k[V]^H \hookrightarrow k[V]$.} the inclusion 
$\CC[V \oplus V]^{\diag^\sigma W} \hookrightarrow \CC[V \oplus V]$, so that
$$
\begin{aligned}
\left( A_{W \times W^\sigma} \right)^{\diag^\sigma W} &:=
  \left(  \CC[V \oplus V] / ( \CC[V \oplus V]^{W \times W^\sigma}_+ )  \right)^{\diag^\sigma W} \\
    & \cong \CC[V \oplus V]^{\diag^\sigma W} / ( \CC[V \oplus V]^{W \times W^\sigma}_+ ).
\end{aligned}
$$

One can then compute
$$
\begin{aligned}
\left\langle
   \left( A_{W \times W^\sigma} \right)_{i,j}, \mathbf{1}
\right\rangle_{\diag^\sigma W}
&= \left\langle 
      \Res^{W \times W^\sigma}_{\diag^\sigma W} 
      \left( \chi_{A_i} \otimes \chi_{A^\sigma_j} \right),   \mathbf{1}
\right\rangle_{\diag^\sigma W} \\
&= \frac{1}{|W|}
     \sum_{w \in W}
       \chi_{A_i}(w) \cdot \chi_{A^\sigma_j}(\sigma(w)) \cdot 1 \\
&= \frac{1}{|W|}
     \sum_{w \in W}
       \chi_{A_i}(w) \cdot   \sigma\left(\chi_{A_j}(w) \right) \\
&= \left\langle
      \chi_{A_i},\overline{\sigma}\chi_{A_j} 
    \right\rangle_W
\end{aligned}
$$
where the first equality uses \eqref{tensor-decomposition}
and the third equality uses Proposition~\ref{coinvariant-trace-comparison}.

To finish the proof, it remains to show that 
$\langle \chi_{A_i},\overline{\sigma}\chi_{A_j} \rangle_W$ is the coefficient of $t^i q^j$ in
$$
W^\sigma(t,q)  =  \sum_\lambda f^\lambda(t) f^{\overline{\sigma}^{-1}(\lambda)}(q).
$$
But this follows from ~\eqref{fake-degree-definition}:  the
class function $\chi_{A_i}$ is the coefficient of $t^i$ in 
$$
\chi_A = \sum_\lambda f^\lambda(t) \cdot \chi^\lambda,
$$
while the class function $\overline{\sigma}\chi_{A_j}$ is the coefficient of $q^j$ in
$$
\overline{\sigma}\chi_{A} 
  =   \sum_\lambda f^\lambda(q) \cdot \chi^{\overline{\sigma}(\lambda)}
   =  \sum_\lambda f^{\overline{\sigma}^{-1}(\lambda)}(q) \cdot \chi^{\lambda}.
$$
\end{proof}

Theorem~\ref{bimahonian-invariant-defn} immediately gives a ``Molien-type'' formula
for $W^\sigma(t,q)$. Note that
\begin{equation}
\label{polynomial-invariants-hilb}
\begin{aligned}
\Hilb(\CC[V \oplus V]^{W \times W^\sigma}; \, t,q)& = \frac{1}{\prod_{i=1}^n (1-t^{d_i})(1-q^{d_i})}.
\end{aligned}
\end{equation}
Recall that $w\in W$ is an element of the group $GL(V)$, so that $w$ and $\sigma(w)$
can be represented by $n\times n$ matrices, thus expressions such as $\det(1-tw)$ 
or $\det(1-q\sigma(w))$ are well-defined polynomials.

\begin{corollary}
\label{Molien-style-formulae}
For any complex reflection group $W$, and Galois automorphism $\sigma$, one has
$$
W^\sigma(t,q) =\frac{1}{|W|}
              \sum_{w \in W} \frac{\prod_{i=1}^n (1-t^{d_i})(1-q^{d_i})}{\det(1-tw)\det(1-q\sigma(w))}.
$$
\end{corollary}
\begin{proof}
Molien's Theorem \cite[\S 2]{Stanley-invariants} tells us that
\begin{equation}
\label{Molien-formula}
 \Hilb(\CC[V \oplus V]^{\diag^\sigma W}; t,q) =
\frac{1}{|W|} \sum_{w \in W} \frac{1}{\det(1-tw)\det(1-q\sigma(w))}.
\end{equation}
On the other hand, $\CC[V \oplus V]^{\diag^\sigma W}$ is a Cohen-Macaulay ring \cite[\S 3]{Stanley-invariants},
and hence a free module over the polynomial subalgebra $\CC[V \oplus V]^{W \times W^\sigma}$.
Thus
\begin{equation}
\label{quotient-of-hilbs}
\begin{aligned}
W^\sigma(t,q) &= \Hilb(\CC[V \oplus V]^{\diag^\sigma W}/
                       (\CC[V \oplus V]^{W \times W^\sigma}_+); \, t,q) \\
       &=  \frac{ \Hilb(\CC[V \oplus V]^{\diag^\sigma W}; \, t,q) }
               { \Hilb(\CC[V \oplus V]^{W \times W^\sigma}; \, t,q) }.
\end{aligned}
\end{equation}
The desired formula for $W^\sigma(t,q)$ results from taking the
quotient of the right side of \eqref{Molien-formula} 
by the right side of \eqref{polynomial-invariants-hilb}.
\end{proof}

It was already noted in \eqref{t-q-symmetry} of the Introduction that
$W^\sigma(q,t)=W^{\sigma^{-1}}(t,q)$, and hence that $W(t,q), \overline{W}(t,q)$
are symmetric polynomials in $t, q$.  It was further noted there that $W^\sigma(1,q)=W^\sigma(q,1)=W(q)$,
so the maximal $q$-degree and maximal $t$-degree in $W^\sigma(t,q)$ are both equal to
the degree of $W(q)$, namely $N^*=\sum_{i=1}^n (d_i-1)$, the number of reflections in $W$.
We note here one further symmetry property enjoyed by $\overline{W}(t,q)$.

\begin{corollary}
\label{symmetry-corollary}
For any complex reflection group $W$,
$$
(tq)^{N^*} \overline{W}(t^{-1},q^{-1}) = \overline{W}(t,q).
$$
In other words, for all $i,j$, 
the monomials $t^i q^j$ and  $t^{N^*-i} q^{N^*-j}$ carry the same coefficient in $\overline{W}(t,q)$.
\end{corollary}
\begin{proof}
We offer two proofs.

\noindent
{\it Proof 1.}
Let $R$ denote the operator on rational functions
$f(t,q)$, defined by 
$$
R(f):=(tq)^{N^*} f(t^{-1},q^{-1}).
$$  
We wish to show that $R$ fixes $\overline{W}(t,q)$. In fact, one
can check that in the formula for $\overline{W}(t,q)$ given by 
Corollary~\ref{Molien-style-formulae}, the operator $R$ will
exchange the summand corresponding to $w$
\begin{equation}
\label{Molien-summand}
\frac{\prod_{i=1}^n (1-t^{d_i})(1-q^{d_i})}
          {\det(1-tw)\det(1-q\overline{w})}
=\frac{\prod_{i=1}^n (1-t^{d_i})(1-q^{d_i})}
          {\det(1-tw)\det(1-qw^{-1})}
\end{equation}
with the summand for $w^{-1}$.  To see this, 
apply $R$ to \eqref{Molien-summand} 
by first sending $t,q \mapsto t^{-1},q^{-1}$,
and then getting rid of all the negative powers of $t, q$ in the numerator and
denominator by multiplying in $(tq)^{N^*}=\frac{(tq)^{\sum_i d_i}}{(tq)^n}$.
The result is 
$$
\frac{\prod_{i=1}^n (t^{d_i}-1)(q^{d_i}-1)}
          {\det(t-w)\det(q-w^{-1})}.
$$
which one can see equals the summand for $w^{-1}$ using these facts:
$$
\begin{aligned}
\det(t-w) &= \det(w) \det(tw^{-1}-1) \\
\det(q-w^{-1}) &= \det(w^{-1})  \det(qw-1) \\
1&=\det(w)\det(w^{-1}).
\end{aligned}
$$

\noindent
{\it Proof 2.}
$\CC[V \oplus V]^{\diag^\sigma W}$ turns out to be a {\it Gorenstein ring},
as we now explain.  In this case, $\diag^\sigma W$ is a subgroup of
$SL(V \oplus V)$, because
$$
\det(w \oplus \overline{w})=\det(w) \overline{\det(w)}=1.
$$
Hence the invariants when $\diag^\sigma W$ acts on $\CC[V \oplus V]$
form a Gorenstein ring by a result of Watanabe; see \cite[Corollary 8.2]{Stanley-invariants}.

This means that the Gorenstein quotient 
$$
\CC[V \oplus V]^{\diag^\sigma W}/(\CC[V \oplus V]^{W \times W^\sigma}_+)
$$
of Krull dimension $0$ is a Poincar\'e duality algebra.  
It shares the same socle bidegree $(N^*,N^*)$
as the larger Gorenstein quotient $\CC[V \oplus V]/(\CC[V \oplus V]^{W \times W^\sigma}_+)$, 
since there is an antidiagonally-invariant element $J(\xx)J(\yy)$ living in this
bidegree.  Here $J(\xx)$ is the {\it Jacobian} 
$$
J(\xx) := \det \left( \frac{\partial f_i}{\partial x_j} \right)_{i,j=1}^n
=\prod_H \left( \ell_H \right)^{s_H},
$$
where $H$ runs through the reflecting hyperplanes for $W$,
with $\ell_H$ any linear form defining $H$, and $s_H$ is the 
number of reflections that fix $H$ pointwise; see \cite[Proposition 4.7]{Stanley-invariants}.
\end{proof}

\begin{example}
\label{cyclic-group-example}
We analyze here in detail the case where $n=1$.  Here $V=\CC^1$
and $W = <c> = \ZZ/d\ZZ$ is cyclic, acting by a representation $\rho: W \rightarrow GL(V)$
that sends $c$ to some primitive $d^{th}$ root-of-unity $\omega$.
Then $W$ acts on $\CC[V]=\CC[x]$ by $c(x)=\omega^{-1}x$,
and $\CC[V]^W=\CC[x^d]$ so there is one fundamental degree $d_1=d$.

The coinvariant algebra is
$$
A=\CC[V]/(\CC[V]^W_+) \cong \CC[x]/(x^d) = \CC\{1,x,x^2,\ldots,x^{d-1}\}
$$
so its Hilbert series gives the Mahonian distribution
$$
W(q)= \Hilb( \CC[x]/(x^d); q) = 1+q+q^2+\cdots +q^{d-1}.
$$
Indexing the irreducible representations/characters as the powers $\{\rho^i\}_{ i \in \ZZ/d\ZZ }$
of the primitive representation $\rho$, one can see that for $i=0,1,\ldots,d-1$
the homogeneous component $A_i=\CC\{x^i\}$ carries the representation $\rho^{d-i}$,
and hence the fake degree polynomials are given by $f^{\mathbf{1}}=1$ and
$f^{\rho^i}(q)=q^{d-i}$ for $i=1,\dots,d-1$.

Note that $W$ is defined over the cyclotomic extension $\QQ[\omega]=\QQ[e^{\frac{2 \pi i}{d}}]$,
for which the Galois automorphisms take the form $\sigma(\omega)=\omega^s$
for some $s \in (\ZZ/d\ZZ)^\times$.  

Working in $\CC[V \oplus V]=\CC[x,y]$, regardless of the choice of $\sigma$,
one has 
$$
\CC[V \oplus V]^{W \times W^\sigma} =\CC[x^d,y^d].
$$
The diagonal invariants $\CC[x,y]^{\diag^\sigma W}$ have
$\CC$-basis given by the monomials 
$$
\{ x^a y^b : a,b \in \NN \text{ and } a+sb \equiv 0 \mod d\}.
$$
Furthermore, $\CC[x,y]^{\diag^\sigma W}$ is a free $\CC[x^d,y^d]$-module with $\CC[x^d,y^d]$-basis
given by the $d$ monomials
$$
\{ x^a y^b : a,b \in \{0,1,\ldots,d-1\}\text{ and } a+sb \equiv 0 \mod d\}.
$$
These monomials form a basis for the quotient 
$$
\CC[V \oplus V]^{\diag^\sigma W}/(\CC[V \oplus V]^{W \times W^\sigma}_+),
$$
and hence
$$
W^\sigma(t,q) = \sum_{\substack{a,b \in \{0,1,\ldots,d-1\}\\ a+sb \equiv 0 \mod n}} t^a q^b.
$$
In the two most important special cases, one has explicitly
$$
\begin{aligned}
\overline{W}(t,q) &= 1+tq+t^2q^2+\cdots+t^{d-1}q^{d-1} \\
       &= \sum_{\lambda} f^\lambda(t) f^\lambda(q) \\
W(t,q)&= 1+t^{d-1}q+t^{d-2}q^2+\cdots+t^2 q^{d-2}+t q^{d-1} \\
       &= \sum_{\lambda} f^\lambda(t) f^{\overline{\lambda}}(q). \\
\end{aligned}
$$
Note that within this family $W=\ZZ/d\ZZ$, one has 
$\overline{W}(t,q)=W(t,q)$ if and only if $d \leq 2$, that is, exactly
when $W$ is actually a {\it real} reflection group-- either $W=\ZZ/2\ZZ$ or the trivial group $W=\ZZ/1\ZZ$.
\end{example}

\section{Bicyclic sieving phenomena}
\label{biCSP-defn-section}

The {\it cyclic sieving phenomenon} for a triple $(X,X(q),C)$ was defined in \cite{RSW}.
Here $X$ is a finite set with an action of a cyclic group $C$, and
$X(q)$ is a polynomial in $q$ with integer coefficients.  We wish to
generalize this notion to actions of bicyclic groups $C \times C'$
and bivariate polynomials $X(t,q)$, so we define this carefully here.

Let $X$ be a finite set with a permutation action of a finite {\it bicyclic group}, that is
a product $C \times C' \cong \ZZ/k\ZZ \times \ZZ/\ell\ZZ$.  Fix embeddings 
$$
\begin{aligned}
\omega: C \hookrightarrow \CC^\times \\
\omega': C' \hookrightarrow \CC^\times
\end{aligned}
$$
into the complex roots-of-unity.  Assume we are given a bivariate polynomial
$X(t,q) \in \ZZ[t,q]$, with nonnegative integer coefficients.

\begin{proposition} (cf. \cite[Proposition 2.1]{RSW})
\label{biCSP-defn-prop}
In the above situation, the following two conditions on the triple $(X,X(t,q),C\times C')$ are equivalent:
\begin{enumerate}
\item[(i)] For any $(c,c') \in C \times C'$ 
$$
X(\omega(c),\omega'(c')) = |\{x \in X: (c,c')x=x\}|.
$$
\item[(ii)] The coefficients $a(i,j)$ uniquely defined by the expansion
$$
X(t,q) \equiv \sum_{\substack{0 \leq i < k\\0 \leq j < \ell}} a(i,j) t^i q^j 
         \mod (t^k-1, q^\ell-1)
$$
have the following interpretation:  $a(i,j)$ is the number of orbits of $C \times C'$ 
on $X$ for which the $C \times C'$-stabilizer subgroup of any element in the orbit lies in the
kernel of the $C \times C'$-character $\rho^{(i,j)}$ defined by
$$
\rho^{(i,j)}(c,c')=\omega(c)^i \omega'(c')^j.
$$
\end{enumerate}
\end{proposition}

\begin{definition}
Say that the triple $(X,X(t,q),C \times C')$ {\it exhibits the bicyclic sieving phenomenon} 
(or {\it biCSP} for short) if it satisfies the conditions of Proposition~\ref{biCSP-defn-prop}.
\end{definition}

In other words, whenever $(X,X(t,q),C \times C')$ exhibits the biCSP, not only is the evaluation 
$X(1,1)$ of $X(t,q)$ telling us the cardinality $|X|$,
but the other root-of-unity evaluations $X(\omega,\omega')$ are telling us all the
$C \times C'$-character values for the permutation action of $C \times C'$ on $X$,
thus determining the representation up to isomorphism.  
Furthermore, the reduction of $X(q,t) \mod (q^k-1,t^\ell-1)$ has a simple interpretation
for its coefficients.

\begin{proof}(of Proposition~\ref{biCSP-defn-prop})
We follow the proof of \cite[Proposition 2.1]{RSW}, by introducing a third equivalent
condition.  

Given $X(t,q)$ define an $\NN^2$-graded complex vector space $A_X = \oplus_{i,j \geq 0} (A_X)_{ij}$
with $\dim_\CC (A_X)_{ij}$ equal to the coefficient of $t^i q^j$ in $X(t,q)$, so that by definition
one has 
$$
\Hilb(A_X; \, t,q)=X(t,q).
$$
Make $C \times C'$ act on $A_X$ by having $(c,c')$ act on the homogeneous component $(A_X)_{ij}$
via the scalar $\rho^{(i,j)}(c,c')=\omega(c)^i \omega'(c')^j$.  

We claim that conditions (i) and (ii) in Proposition~\ref{biCSP-defn-prop} are equivalent to
\begin{enumerate}
\item[(iii)] $A_X \cong \CC[X]$ as (ungraded) $C \times C'$-representations.
\end{enumerate}

The equivalence of (i) and (iii) follows immediately from the observation that
$C \times C'$-representations are determined by their character values for each
element $(c,c')$;  in $A_X$ this character value is $X(\omega(c),\omega'(c'))$ by
construction, and in $\CC[X]$ this character value is $|\{x \in X: (c,c')x=x\}|$.

For the equivalence of (ii) and (iii), first note that the complete set of irreducible
representations or characters of $C \times C'$ are given by 
$\{\rho^{(i,j)}\}$ for $0 \leq i < k$ and $0 \leq j< \ell$.
Consequently, (ii) holds if and only if for all such $i,j$ one has
\begin{equation}
\label{biCSP-multiplicity-equality}
\langle \rho^{(i,j)} , A_X \rangle_{C \times C'} =
\langle \rho^{(i,j)} , \CC[X] \rangle_{C \times C'}.
\end{equation}
We compute the left side of \eqref{biCSP-multiplicity-equality}:
$$
\begin{aligned}
\langle \rho^{(i,j)} , A_X \rangle_{C \times C'}
&= \frac{1}{|C \times C'|} \sum_{(c,c') \in C \times C'} \rho^{(i,j)}((c,c')^{-1}) \chi_{A_X}(c,c') \\
&= \frac{1}{k \ell} \sum_{\substack{(\omega,\omega'):\\ \omega^k=1\\ (\omega')^\ell=1}} 
         \omega^{-i} (\omega')^{-j} X(\omega,\omega') \\
&= a(i,j).
\end{aligned}
$$
To compute the right side of \eqref{biCSP-multiplicity-equality},
first decompose $X$ into its various $C \times C'$-orbits $\OOO$.
Denote by $G_\OOO$ the $C \times C'$-stabilizer subgroup of any element in
the orbit $\OOO$.  One then has
$$
\begin{aligned}
\langle \rho^{(i,j)} , \CC[X] \rangle_{C \times C'}
&= \sum_{ \OOO } \langle \rho^{(i,j)} , \CC[\OOO] \rangle_{C \times C'} \\
&= \sum_{ \OOO } \langle \rho^{(i,j)} , \Ind_{G_\OOO}^{C \times C'} {\mathbf 1} \rangle_{C \times C'} \\
&= \sum_{ \OOO } \langle \Res^{C \times C'}_{G_\OOO} \rho^{(i,j)} , {\mathbf 1} \rangle_{G_\OOO} \\
\end{aligned}
$$
and each term in this last sum is either $1$ or $0$ depending upon whether
$G_\OOO$ lies in the kernel of $\rho^{(i,j)}$, or not.  Thus (ii) holds if and only
if \eqref{biCSP-multiplicity-equality} holds for all $i,j$, 
that is, if and only if (iii) holds.
\end{proof}

\section{Springer's regular elements and proof of Theorem~\ref{biCSP-theorem}}
\label{Springer-section}

In order to prove Theorem~\ref{biCSP-theorem}, we first
must recall Springer's notion of a regular element in a complex reflection group.

\begin{definition}
Given a complex reflection group $W$ in $GL(V)$, a vector $v \in V$
is called {\it regular} if it lies on none of the reflecting hyperplanes
for reflections in $W$.  An element $c$ in $W$ is called {\it regular} if
it has a regular eigenvector $v$; the eigenvalue $\omega$ for
which $c(v)=\omega v$ will be called the accompanying {\it regular} eigenvalue.
One can show \cite[p. 170]{Springer} that $\omega$ will have the same
multiplicative order in $\CC^\times$ as the multiplicative order of $c$ in $W$.
\end{definition}

Note that the eigenvalues of a regular element need not lie in the {\it smallest}
cyclotomic extension over which $W$ is defined;  Example~\ref{type-A-regular-elements} 
below illustrates this.  However, if one chooses $m$ to be the least common multiple
of the orders of the elements of $W$, then any such eigenvalue {\it will} lie in the
cyclotomic extension $\QQ[e^{\frac{2\pi i}{m}}]$.  Hence any $\sigma$ in 
$\Gal(\QQ[e^{\frac{2\pi i}{m}}]/\QQ)$ acts on any such eigenvalue $\omega$, with
$\sigma(\omega)=\omega^s$ for some unique $s \in (\ZZ/d\ZZ)^\times$, where $d$ is the order
of the regular element $c$.  For the remainder
of this section we will always assume that $m$ is this least common multiple.

\begin{example}
\label{type-A-regular-elements}
Springer \cite[\S 5]{Springer} classified the regular elements
in the real reflection groups.  When $W=\symm_n$, regular vectors are
those vectors in $\CC^n$ for which all coordinates are distinct.
Any $n$-cycle or $(n-1)$-cycle is a regular element--
for example $c=(1 \, 2 \, \cdots \, n)$ has regular eigenvalue $\omega$
and regular eigenvector $v=(1,\omega,\omega^2,\ldots,\omega^{n-1})$, if 
$\omega$ is any 
primitive $n^{th}$ root-of-unity.  
Similarly, $c'=(1 \, 2 \, \cdots \, n-1)(n)$ has
regular eigenvalue $\zeta$ and regular eigenvector $v'=(1,\zeta,\zeta^2,\ldots,\zeta^{n-1},0)$ for any
primitive $(n-1)^{st}$ root-of-unity $\zeta$.  Any power of an $n$-cycle will
therefore also be regular, with the same regular eigenvector $v$, as will
any power of an $(n-1)$-cycle.  

It turns out (and is not hard to see directly)
that there are no other regular elements besides powers of $n$-cycles and
$(n-1)$-cycles in $\symm_n$.
\end{example}

Springer proved the following.

\begin{theorem}\cite[Prop. 4.5]{Springer}
\label{Springer-theorem}
For a regular element $c$ in a complex reflection group $W$, with
an accompanying regular eigenvalue $\omega$, one has
$$
\chi^\lambda(c) = f^\lambda(\omega^{-1}).
$$
\end{theorem}

This leads to an interesting interaction between Galois conjugates $\sigma$
and a regular element $c$.  Recall that there is a unique $s \in (\ZZ/d\ZZ)^\times$
such that the regular eigenvalue $\omega$ satisfies $\sigma(\omega)=\omega^s$.

\begin{proposition}
\label{sigma-power}
With the above notation, for any $W$-irreducible $\lambda$, one has
$$
\chi^{\sigma(\lambda)}(c) = \chi^\lambda(c^s).
$$
\end{proposition}  
\begin{proof}
Note that since $c$ is regular with regular eigenvalue $\omega$,
one has that $c^s$ is regular with regular eigenvalue $\omega^s$.
Then one has
$$
\begin{aligned}
\chi^{\sigma(\lambda)}(c) 
   &= \sigma( \chi^\lambda(c) ) \\
   &= \sigma \left( f^\lambda( \omega^{-1} ) \right) \\
   &= f^\lambda \left( \sigma(\omega^{-1} ) \right) \\
   &= f^\lambda( \omega^{-s} ) \\
   &= \chi^{\lambda}( c^s )
\end{aligned}
$$
where the second and fifth equalities use Theorem~\ref{Springer-theorem}.
\end{proof}

From this we can now prove Theorem~\ref{biCSP-theorem}, which we rephrase here as
a biCSP, after establishing the appropriate notation.  

Given two regular elements $c,c'$ in $W$,
with regular eigenvalues $\omega, \omega'$, embed the cyclic groups
$C=\langle c \rangle, C'=\langle c' \rangle$ into $\CC^\times$ by
sending 
$$
\begin{aligned}
c &\mapsto \omega^{-1} \\
c' &\mapsto (\omega')^{-1}.
\end{aligned}
$$

As above, let $d$ denote the multiplicative order of $c$ and of $\omega$,
and given the Galois automorphism $\sigma$, define $s \in (\ZZ/d\ZZ)^\times$
by the property that $\sigma(\omega)=\omega^s$

Let $X=W$, carrying the
following $\sigma$-twisted left-action of $C \times C'$:
$$
(c,c') \cdot w := c^s w (c')^{-1}.
$$
Let $X(t,q)=W^\sigma(t,q)$.

\vskip .2in
\noindent
{\bf Theorem~\ref{biCSP-theorem} (rephrased)}.
{\it
In the above setting, $(X,X(t,q),C \times C')$ exhibits the biCSP.
}
\begin{proof}
As noted earlier, if $c$ is a regular element $c$ of $W$ with regular eigenvalue $\omega$,
then any power $c^i$ is also regular, with regular eigenvalue $\omega^i$.  Hence
it suffices for us to show that
$$
W^\sigma(\omega^{-1},(\omega')^{-1})=|\{w \in W: c^s w (c')^{-1}=w\}|.
$$
This follows from a string of equalities, explained below:
$$
\begin{aligned}
W^\sigma(\omega^{-1},(\omega')^{-1}) 
   &= \sum_\lambda f^{\sigma(\lambda)}(\omega^{-1}) f^{\overline{\lambda}}((\omega')^{-1})\\
   &= \sum_\lambda \chi^{\sigma(\lambda)}(c) \chi^{\overline{\lambda}}(c')\\
   &= \sum_\lambda \chi^\lambda(c^s) \overline{\chi^\lambda(c')}\\
   &=\begin{cases}
          |\Cent_{W}(c^s)| &\text{ if } c', c^s, \text{ are }W\text{-conjugate} \\
          0                   &\text{ otherwise.}\\
                    \end{cases} \\
&= |\{w \in W: w^{-1} c^s w = c' \}| \\
&= |\{w \in W: c^s w (c')^{-1} = w \}|
\end{aligned}
$$
The first equality is by Definition~\ref{bimahonian-fake-degree-defn}
for $W^\sigma(t,q)$.  
The second equality uses Theorem~\ref{Springer-theorem}.
The third equality uses Proposition~\ref{sigma-power}.
The fourth equality uses the column orthogonality relation for the character table
of $W$.
\end{proof}

%
%
%

\section{Type A:  bipartite partitions and work of Carlitz, Wright, Gordon} 
\label{Gordon-section}

We discuss briefly here the known root-of-unity evaluations for $W^\sigma(t,q)$
in the much-studied case $W=\symm_n$.  Since this $W$ is defined over $\QQ$, one
may assume $\sigma$ is the identity and study only $\symm_n(t,q)=W(t,q)$.

Theorem~\ref{biCSP-theorem} combinatorially interprets $\symm_n(\omega,\omega')$ when
the roots of unity in question have orders that divide either $n-1$ or $n$, since these are exactly
the orders of regular elements in $W=\symm_n$;  see Example~\ref{type-A-regular-elements}.  
This combinatorial interpretation
is consistent with evaluations done by Gordon \cite{Gordon}, building on work of Carlitz \cite{Carlitz}
and Wright \cite{Wright}.  In fact, Gordon's work actually leads to some further evaluations of $\symm_n(\omega,\omega')$
more generally when the orders of $\omega, \omega'$ are at most $n$, as we now explain.

The starting point both for Carlitz and Wright is a generating function for
bipartite partitions.  Every $\diag \symm_n$-orbit of 
monomials $\xx^{\mathbf i} \yy^{\mathbf j}=\prod_{m=1}^n x_m^{i_m} y_m^{j_m}$ in 
$
\CC[V \oplus V]=\CC[x_1,\ldots,x_n,y_1,\ldots,y_n]
$
corresponds to an (unordered) multiset of $n$ exponent vectors $\{(i_m,j_m)\}_{m=1,\ldots,n}$.
A canonically chosen ordering for this multiset, say in lexicographic order,
is known as a {\it bipartite partition}; 
for further background on this topic, see Roselle \cite{Ros}, 
Solomon \cite{Solomon}, Garsia and Gessel \cite{GarsiaGessel}, 
Stanley \cite[Example 5.3]{Stanley-invariants} and 
Adin-Gessel-Roichman \cite[\S4]{AGR}. 

Since such orbits of monomials also form a $\CC$-basis for $\CC[V \oplus V]^{\diag \symm_n}$, one obtains
the following generating function for the $\symm_n(t,q)$:

\begin{equation}
\label{starting-gf}
\begin{aligned}
\prod_{i,j \geq 0} \frac{1}{1-t^i q^j u}
   &= \sum_{n \geq 0} u^n \Hilb( \CC[V \oplus V]^{\diag \symm_n} ; \, t,q) \\
   &= \sum_{n \geq 0} u^n \frac{\symm_n(t,q)}{(t;t)_n (q;q)_n} 
\end{aligned}
\end{equation}
where we have used the notation
$$
(q;q)_n := (1-q)(1-q^2) \cdots (1-q^n) = \frac{1}{\Hilb(\CC[V]^{\symm_n}, q)}
$$
to rewrite \eqref{quotient-of-hilbs} in this situation.
Logarithmically differentiating \eqref{starting-gf} gives the following recurrence 
used by Wright \cite{Wright}, Carlitz \cite{Carlitz}, and then Gordon \cite{Gordon}.

\begin{proposition}
\label{Wright-recursion}
$\symm_n(t,q)$ is uniquely defined by the recurrence
$$
\symm_n(t,q) = \frac{1}{n} \sum_{m=1}^n 
             \frac{(t;t)_n (q;q)_n}{(t;t)_{n-m} (q;q)_{n-m}}
                     \frac{\symm_{n-m}(t,q)}{(1-t^m)(1-q^m)}
$$
with initial condition $\symm_0(t,q)=1$. 
\end{proposition}

From Proposition~\ref{Wright-recursion}, Gordon deduced via induction on $n$
the following lemma about the behavior of $\symm_n(t,q)$ when one of the two
variables is specialized to a root of unity.

\begin{lemma}\cite[\S3]{Gordon}
\label{Gordon-lemma}
Let $\omega$ be a primitive $\ell^{th}$ root of unity, and write $n =m \ell  + r$
with $0 \leq r < \ell$.  Then
$$
\symm_n(\omega,q) = \frac{(q;q)_n}{(q;q)_r (1-q^{\ell})^{m}} \symm_r(\omega,q).
$$
\end{lemma}

The next corollary gives two interesting evaluations of $\symm_n(\omega,\omega')$
obtainable with a little work by taking $q=\omega'$ in Lemma~\ref{Gordon-lemma}.
The first was already deduced by Gordon as one of his main results.  Although the 
second he seems not to have written down, we omit its relatively straightforward proof here.

\begin{corollary}
\label{Gordon-evaluations-corollary}
Let $\omega, \omega'$ be roots of unity.
\begin{enumerate}
\item[(i)]
If they have unequal orders, with both orders at most $n$, then 
$$
\symm_n(\omega, \omega')=0.
$$
\item[(ii)]
If their orders are both $\ell$, then
$$
\symm_n(\omega,\omega') = \ell^m m! \symm_r(\omega,\omega')
$$
where one has written uniquely $n =m \ell  + r$ with $0 \leq r < \ell$.

In particular, if $n \equiv 0,1 \mod \ell$ then $\symm_n(\omega,\omega') = \ell^m m!.$
\end{enumerate}
\end{corollary}

We explain here how Corollary~\ref{Gordon-evaluations-corollary}
recovers (more than) the assertion of Theorem~\ref{biCSP-theorem} for $W=\symm_n$.
Assume that $c, c'$ in $W=\symm_n$ are regular elements,  and $\omega, \omega'$ their 
corresponding regular eigenvalues.  According to Example~\ref{type-A-regular-elements}, 
$c, c'$ must each be a power
of an $(n-1)$-cycle or $n$-cycle.  Note that this means that
$c, c'$ will be $W$-conjugate (that is, have the same cycle type) 
if and only if they have the same multiplicative order.  

\vskip .1in
\noindent
{\bf Case 1.}  $c,c'$ are not $W$-conjugate.

In this case,  Theorem~\ref{biCSP-theorem} predicts $\symm_n(\omega,\omega')=0$.
By the above comment, $\omega, \omega'$ have unequal orders, and hence
Corollary~\ref{Gordon-evaluations-corollary}(i) 
also predicts $\symm_n(\omega, \omega')=0$.

\vskip .1in
\noindent
{\bf Case 2.}  $c,c'$ are $W$-conjugate.

In this case,  Theorem~\ref{biCSP-theorem}
predicts that $\symm_n(\omega',\omega)$ will be the number of elements
$w \in W$ with $cw(c')^{-1}=w$, or equivalently, $w^{-1}cw=c'$,  
which is counted by the centralizer-order $|\Cent_W(c)|=\ell^m m!$.
If $c, c'$ have order $\ell$, their being regular forces 
$n=m\ell+r$ with remainder $r=0$ or $r=1$, and 
Corollary~\ref{Gordon-evaluations-corollary}(ii) predicts
$$
\symm_n(\omega',\omega) = \ell^m m! \symm_r(\omega',\omega)= \ell^m m!
$$
since $\symm_0(t,q)=\symm_1(t,q)=1$.

We close this section with a few remarks.

\begin{remark} \rm \
The results in \cite[\S 4, 5]{BSS} are assertions about
the coefficients $a_{d,d}(i,j)$ defined in \eqref{reduced-coefficients-defn}
for $W=\symm_n$ with $n \cong 0,1 \mod d$.  The approach taken there is
to deduce them from an explicit number-theoretic formula \cite[Theorem 4.1]{BSS} for
$a_{n,n}(i,j)$.  This formula is in turn derived using a result of Stanley and 
Kraskiewicz-Weyman (equivalent to Theorem~\ref{Springer-theorem} in type $A$)
by first reinterpreting the coefficients $a_{n,n}(i,j)$ as certain intertwining numbers:
\begin{equation}
\label{BSS-inner-product}
a_{n,n}(i,j)=
\langle 
\Ind_{C_n}^{\symm_n }\rho^{i}, \Ind_{C_n}^{\symm_n} \rho^{j} 
\rangle_{\symm_n}.
\end{equation}
Here $C_n$ denotes the cyclic subgroup of $\symm_n$ generated by an $n$-cycle $c$,
and $\rho: C_n \rightarrow \CC^\times$ is the primitive character sending $c$ 
to a primitive $n^{th}$ root-of-unity $\omega$.

It is not hard to check that the combinatorial interpretation for $a_{n,n}(i,j)$ given by
Theorem~\ref{biCSP-theorem} (counting $C_n \times C_n$-orbits $\OOO$ on $\symm_n$ 
whose stabilizer $G_\OOO \subset \ker \rho^{(i,j)}$)
is equivalent to \eqref{BSS-inner-product} via Mackey's formula.  Thus the combinatorial
interpretation for $a_{n,n}(i,j)$ gives another 
route to the results of \cite[\S 4,5]{BSS}.
\end{remark}

\begin{remark}
The identity \eqref{starting-gf} gives a concise generating function that
compiles bivariate distributions $(\maj(w), \maj(w^{-1}))$ for all of the
symmetric groups $W=\symm_n$.  It is a specialization of a stronger identity
due to Garsia and Gessel \cite{GarsiaGessel} that does the same for the four-variate
distribution of $(\maj(w), \maj(w^{-1}), \des(w), \des(w^{-1}))$ where
$\des(w):=\{i: w(i) > w(i+1)\}$ is the number of {\it descents} in $w$.

More recently, Foata and Han \cite[Eqn. (1.8)]{FoataHan} generalized this 
by giving such a generating function for the {\it hyperoctahedral groups} 
$W=W(B_n)=\ZZ/2\ZZ \wr \symm_n$ of signed permutations.  Their result
incorporates the five-variate distribution of certain statistics
$(\fmaj(w), \fmaj(w^{-1}), \fdes(w), \fdes(w^{-1}), \nneg(w))$.  Here $\nneg(w)$ is the
number of negative signs appearing in the signed permutation, so this distribution
for $W(B_n)$ can be specialized to the previous one for $\symm_n$.
\end{remark}

\begin{remark} \rm \ 
Corollary~\ref{Gordon-evaluations-corollary} 
leaves us with the question of what can one say about 
$\symm_r(\omega,\omega')$ for roots of unity $\omega, \omega'$ of the {\it same} order $d$
when $2 \leq r \leq d-1$.  In general, 
such evaluations $\symm_r(\omega,\omega')$ 
can be negative real numbers, or can lie in $\CC \setminus \RR$.

There is however at least one more (somewhat trivial) thing one can say about
the ``antidiagonal'' values.

\begin{proposition}
For any root of unity $\omega$, one has that $\symm_n(\omega,\omega^{-1})$ lies
in $\RR$.
\end{proposition}
\begin{proof}
The fact that $\symm_n(t,q)$ has real coefficients and is symmetric in $t,q$ 
implies that $\symm_n(\omega,\omega^{-1})$
is fixed under complex conjugation:
$$
\overline{\symm_n(\omega,\omega^{-1})} 
  = \symm_n(\overline{\omega},\overline{\omega^{-1}})
    = \symm_n(\omega^{-1},\omega)
     = \symm_n(\omega, \omega^{-1}).
$$
\end{proof}
\end{remark}

%

\section{The wreath products $\ZZ/d\ZZ \wr \symm_n$, tableaux, and the flag-major index}
\label{fmaj-section}

  We review here for the classical complex reflection groups $W=\ZZ/d\ZZ \wr \symm_n$, how
one expresses the fake degrees and bimahonian distributions in terms of major-index-like
statistics, both on tableaux and on $W$.  The one new result here is Theorem~\ref{bimahonian-by-fmaj-theorem},
which generalizes a result of Adin and Roichman \cite{AdinRoichman}.

   We first review the motivating special case when $W=\symm_n$;  see also
\cite{FoataSchutzenberger}.  Recall that the irreducible complex representations of 
$\symm_n$ are indexed by partitions $\lambda$ of $n$. We will use 
$\lambda$ to denote both the irreducible representation and the partition. 
It appears that Lusztig first computed (see
\cite[Prop. 4.11]{Stanley-invariants}, and \cite{KraskiewiczWeyman}) 
the following formula for the fake degree polynomial $f^\lambda(q)$ for
the representation $\lambda$:
$$
f^\lambda(q) =\sum_{Q \in \SYT(\lambda)} q^{\maj(Q)}.
$$
Here $\SYT(\lambda)$ is the set of {\it standard Young
tableaux} of shape $\lambda$, that is, fillings of the Ferrers/Young diagram
for $\lambda$ with each number $1,2,\ldots,n$ occurring exactly once,
increasing left-to-right in rows and top-to-bottom in columns.
The {\it major index} statistic $\maj(Q)$ is defined by
$$
\maj(Q):=\sum_{ i \in \Des(Q)} i,
$$
where the {\it descent set} $\Des(Q)$ is the set of values $i$ for which $i+1$ occurs
in a lower row of $Q$.
The relation to the bimahonian distribution is provided by the {\it Robinson-Schensted
correspondence}, which gives a bijection between permutations 
$w$ in $\symm_n$ and pairs $(P,Q)$ of
standard Young tableaux of the same shape.  Some fundamental properties of this
bijection are that if $w \mapsto (P,Q)$, then 
\begin{equation}
\label{RSK-properties}
\begin{aligned}
w^{-1} &\mapsto (Q,P) \\
\Des(w) &= \Des(Q) \\
\Des(w^{-1})& =\Des(P).
\end{aligned}
\end{equation}

 Hence one obtains
\begin{equation}
\label{type-A-RSK}
\begin{aligned}
W(t,q) &= \sum_\lambda f^\lambda(t) f^\lambda(q)  \\
       &= \sum_{(P,Q)} t^{\maj(Q)} q^{\maj(P)} \\
       &= \sum_{w \in \symm_n} t^{\maj(w)} q^{\maj(w^{-1})} 
\end{aligned}
\end{equation}
and where the second sum is over all pairs $(P,Q)$ of standard Young tableaux
of size $n$ of the same shape.

  For $W$ a Weyl group of type $B_n$, that is, $W=\ZZ/2\ZZ \wr \symm_n$,
Adin and Roichman \cite{AdinRoichman}, defined a {\it flag major index} statistic $\fmaj(w)$, 
equidistributed with the Coxeter group length.  More generally, they define such a statistic
for $W=\ZZ/d\ZZ \wr \symm_n$ with the property that
\begin{equation}
\label{Adin-Roichman-fmaj}
\begin{aligned}
 \Hilb(\CC[V \oplus V]^{\Delta W}/ \CC[V \oplus V]^{W \times W}; q,t) &\left( = W(t,q) \right) \\
 &  =\sum_{w \in \ZZ/d\ZZ \wr \symm_n} t^{\fmaj(w)} q^{\fmaj(w^{-1})}.
\end{aligned}
\end{equation}
A variation on the $d=2$ (type $B$) case of this appears in
Biagioli and Bergeron \cite{BergeronBiagioli}.
For $W$ a Weyl group of type $D_n$, a similar $\fmaj(w)$ statistic was defined by
Biagioli and Caselli \cite{BiagioliCaselli}, who proved an analogous
result to \eqref{Adin-Roichman-fmaj}.

As one might expect, there is an approach to \eqref{Adin-Roichman-fmaj} as in
\eqref{type-A-RSK}:  start with the fake-degree Definition~\ref{bimahonian-fake-degree-defn}
for $W(t,q)$, use a tableau formula for the fake-degrees, and then
apply something like a Robinson-Schensted correspondence-- such a
proof is sketched in \cite[\S 5]{AdinRoichman-FPSAC}.  We pursue here a
similar approach\footnote{The second author thanks M. Ta\c{s}kin 
for explaining how this works in type $B_n$.} 
to the more general Theorem~\ref{bimahonian-by-fmaj-theorem}.

Let $\omega$ be a primitive $d^{th}$ root-of-unity.   A typical element
$w$ in $W=\ZZ/d\ZZ \wr \symm_n$, sending $e_i$ to $\omega^k e_j$, can be
represented by a string of 
letters $w_1\cdots w_n$, where $w_i = \omega^k j$.  Here the letters come from an
alphabet that we linearly order as follows:
\begin{equation}
\label{ordered-alphabet}
\begin{aligned}
\omega^{d-1} 1 &\prec \cdots \prec \omega^{d-1} n \\
               &\prec \cdots \\
\prec \omega^{1} 1 &\prec \cdots \prec \omega^{1} n \\
\prec \omega^{0} 1 &\prec \cdots \prec \omega^{0} n.
\end{aligned}
\end{equation}
Call the letters $\omega^k 1,\ldots,\omega^k n$ the
{\it $k^{th}$ subalphabet}, and let $r_k(w)$ denote the number of letters $w_i$ that
lie in this subalphabet. Let
$$
\maj(w):=\sum_{\substack{i=1,\ldots ,n-1:\\ w_{i+1} \prec w_i}} i
$$
denote the usual major index of $w$ with respect to this order, as defined
in \eqref{maj-definition} above.  Although $\fmaj(w)$ can be defined in
a somewhat more algebraic way, it is shown in \cite[Theorem 3.1]{AdinRoichman} that it is
equivalent to this combinatorial expression:
$$
\fmaj(w)= d \cdot \maj(w) + \sum_{k=0}^{d-1} k \cdot r_k(w).
$$

The smallest cyclotomic extension over which $W=\ZZ/d\ZZ \wr \symm_n$ can be defined is 
$\QQ[e^{\frac{2 \pi i}{d}}]$.  
Note that since a Galois automorphism $\sigma \in \Gal(\QQ[e^{\frac{2 \pi i}{d}}]/\QQ)$ will
always be of the form $\sigma(\omega)=\omega^s$ for some $s \in (\ZZ/d\ZZ)^\times$,
one has the interesting feature here that $\sigma(w) \in W$ for all $w \in W$.
Thus $W^\sigma = \sigma(W) = W$, although $\sigma$ does not fix $w$ pointwise.

Tableaux expressions for the fake degrees of Weyl groups of type
$B_n, D_n$ were computed originally by Lusztig \cite{Lusztig}, and
generalized to the Shepard-Todd infinite family $G(de,e,n)$ 
of complex reflection groups by Stembridge \cite{Stembridge}.
Irreducibles for $W= \ZZ/d\ZZ \wr \symm_n$ can be indexed
by skew diagrams 
$$
\lambda=(\lambda^{d-1} \oplus \cdots \oplus \lambda^1 \oplus \lambda^0)
$$
having $n$ cells total, consisting of $d$-tuples of Ferrers diagrams $\lambda^k$, arranged
in the plane so that $\lambda^{k-1}$ is northeast of $\lambda^k$ (using English notation
for Ferrers diagrams).  Define a {\it standard Young tableau} of shape $\lambda$ to be
a filling of the skew diagram $\lambda$ with the numbers $1,2,\ldots,n$, increasing
left-to-right in rows and top-to-bottom in columns.  Given such a standard Young tableau $Q$
of the skew shape $\lambda$ 
define its descent set $\Des(Q)$ and major index $\maj(w)$ with respect to the ordering
\eqref{ordered-alphabet}, and then define
\begin{equation}
\label{fmaj-of-tableau-defn}
\fmaj(Q)= d \cdot \maj(Q) + \sum_{k=0}^{d-1} k \cdot |\lambda^k|.
\end{equation}
After reviewing this indexing of irreducibles, Stembridge proves the following.

\begin{theorem}(\cite[Theorem 5.3]{Stembridge})
\label{Stembridge-fake-degree}
For  $W= \ZZ/d\ZZ \wr \symm_n$, one has the fake degree expression
$$
f^{\lambda} = \sum_{Q} q^{\fmaj (Q)}.
$$
where $Q$ runs over all standard tableaux of shape $\lambda$, 
\end{theorem}

Lastly, we recall from \cite{StantonWhite} one of the standard generalizations of
the Robinson-Schensted correspondence to $W=\ZZ/d\ZZ \wr \symm_n$.
Given $w \in W$, one produces a pair $(P,Q)$ of standard tableaux of the same skew shape 
$\lambda$ having 
\begin{equation}
\label{r-relationship}
|\lambda_k|=r_k(w)
\end{equation}
in which
the letters $w_i$ coming from the $k^{th}$ subalphabet are inserted using the
usual Robinson-Schensted algorithm into the subtableaux
$P^k$ that occupies the subshape $\lambda^k$, and their position $i$ is recorded
in the subtableaux $Q^k$.  As one varies over all $\lambda$
with $n$ cells as above, one gets a bijection to $W$.

\begin{proposition}
\label{d-tuple-RSK-properties}
This bijection $w \mapsto (P,Q)$ between 
$W=\ZZ/d\ZZ \wr \symm_n$ and pairs of standard Young tableaux of the
same shape $\lambda=(\lambda^{d-1} \oplus \cdots \oplus \lambda^1 \oplus \lambda^0)$ having $n$ cells,
has the following properties:
\begin{enumerate}
\item[(i)]
$$
\overline{w}^{-1} \mapsto (Q,P).
$$
\item[(ii)]
$$
\begin{aligned}
\Des(w) &= \Des(Q)\\
\fmaj(Q)&=\fmaj(w) \\
\fmaj(P)&=\fmaj(\overline{w}^{-1}).
\end{aligned}
$$
\item[(iii)] For any Galois automorphism $\sigma$,
$$
\sigma(w) \mapsto (\sigma(P), \sigma(Q)).
$$
\end{enumerate}
Here $\sigma(Q)$ denotes the skew tableaux 
obtained from $Q$ by re-indexing its subtableaux $Q^i$ according
to the permutation by which $\sigma$, thought of as an element of $(\ZZ/d\ZZ)^\times$,
acts on the indices $i \in \ZZ/d\ZZ$.
\end{proposition}

\begin{proof}
The first two properties follow immediately from \eqref{RSK-properties},
\eqref{fmaj-of-tableau-defn}, and \eqref{r-relationship},
while the third is an easy consequence of the definitions.  
\end{proof}

We can now recall the statement of Theorem~\ref{bimahonian-by-fmaj-theorem}
and prove it.
\vskip .1in
\noindent
{\bf Theorem~\ref{bimahonian-by-fmaj-theorem}.}
{\it 
For the wreath products $W=\ZZ/d\ZZ \wr \symm_n$ and
any Galois automorphism $\sigma \in \Gal(\QQ[e^{\frac{2\pi i}{d}}]/\QQ)$, one has
$$
W^\sigma(t,q) = \sum_{w \in W} q^{\fmaj(w)} t^{\fmaj(\sigma(w^{-1}))}.
$$
}
\begin{proof}
One calculates as follows, using the above bijection $w \mapsto (P,Q)$
and its properties from Proposition~\ref{d-tuple-RSK-properties}:
$$
\begin{aligned}
\sum_{w \in W} q^{\fmaj(w)} t^{\fmaj(\sigma(w^{-1})  )} 
  &= \sum_{w \in W} q^{\fmaj(\overline{\sigma}^{-1}(w))} t^{\fmaj( \overline{w}^{-1}  )} \\
  &= \sum_{(P,Q)}  q^{\fmaj(\overline{\sigma}^{-1}(Q))} t^{\fmaj(P)}  \\
  &= \sum_\lambda f^{\overline{\sigma}^{-1}(\lambda)}(q) f^{\lambda}(t) \\
  &= W^\sigma(t,q) 
\end{aligned}
$$
\end{proof}

\begin{question}
For which complex reflection groups can one produce an $\fmaj$ statistic that
allows one to generalize Theorem~\ref{bimahonian-by-fmaj-theorem}?
\end{question}

\noindent
As a first step, one might try to do this for 
the infinite family $G(de,e,n)$ of complex reflection groups.
Such a generalization might involve the fake-degree formulae of Stembridge \cite{Stembridge}
for $G(de,e,n)$.
%

\section*{Acknowledgements}
The authors thank Stephen Griffeth and M\"uge Ta\c{s}kin for helpful conversations.

\end{document}